\def\C{{\cal C}}
\baselineskip=14pt
\parskip=10pt
\def\Tilde{\char126\relax}
\def\halmos{\hbox{\vrule height0.15cm width0.01cm\vbox{\hrule height
 0.01cm width0.2cm \vskip0.15cm \hrule height 0.01cm width0.2cm}\vrule
 height0.15cm width 0.01cm}}
\font\eightrm=cmr8  
\font\eighttt=cmtt8
\magnification=\magstephalf

\parindent=0pt
\overfullrule=0in
\bf
\centerline{AUTOMATED COUNTING of LEGO TOWERS }
\medskip
\it
\centerline{Doron ZEILBERGER \footnote{$^1$}
{\baselineskip=9pt
\eightrm  \raggedright
Department of Mathematics, Temple University,
Philadelphia, PA 19122, USA. 
{\baselineskip=9pt
\eighttt zeilberg@math.temple.edu ; \break http://www.math.temple.edu/\Tilde 
zeilberg. }
This paper
is accompanied by a Maple package, {\eighttt LEGO}, downloadable
from the author's website. LEGO is a registered trademark.
Supported in part by the NSF. Jan. 31, 1997.  This version: 10/14/97.
}}
\medskip
\qquad\qquad
{Dedicated to Henry Wadsworth Gould on his turning $F_{10}+F_7+F_3$ years
young}
\medskip
\rm
{\bf Abstract:}
H.N.V. Temperley's method for counting vertically
convex polyominoes is modified, generalized, and most importantly,
programmed (in Maple).
\medskip
{\bf Preface}
 
I have never met Henry Gould in person, but have always admired
his work. In particular, his charming article
[G2] served me well many times, his (so far non-shaloshable) famous
identities awed and frustrated me, and I still hope to use
his work on Euler's constant[G3] to, who knows?, prove that it is
irrational.
I hope that this paper will please Henry, especially since his
favorite sequence[G1] makes a brief (`cameo role') appearance.
 
{\bf Toys and Toy Models}
 
In spite of the many great triumphs of mathematics and science, there are
many more problems that we {\it can't} solve than ones that we 
{\it can}. One of
those, that so far defied us,
is that of enumerating {\it animals}.
Even Viennot's[V] powerful theory of heaps, that was so successful
in enumerating  {\it directed animals} (with the deceptively
simple formula $3^n$), seems, at present, to be incapable of counting
plain animals.
 
A two-dimensional {\it animal}, alias {\it polyomino}, can be realized
in terms of a {\it LEGO tower}.
Suppose that we have an infinite supply of $1 \times a$ ($a \geq 1$)
LEGO pieces. Then every floor of the tower consists of a finite
horizontal sequence of pieces separated by gaps. A 
vertical sequence of floors
constitutes an animal if the resulting configuration is {\it connected}.
Each floor can be described by the sequence of lengths of the pieces
intertwined by the sequence of lengths of the gaps, so we have
an infinite alphabet, each letter being of the form:
$a_1, b_1 , a_2, b_2, \dots , a_n$. 
Fixing the leftmost square of the bottom floor at the origin,
one can view an animal as a word in this infinite alphabet, together
with a specification of `interfaces', which indicates where to place
the leftmost square of the next floor in relation to the floor below it.
The resulting creature is an animal iff the resulting configuration is
connected. Note that in addition to the complication of having
an infinite ( in fact $\sum_{k=0}^{\infty} \infty^{2k+1}$) 
alphabet, the condition of
connectedness is global which makes it strongly non-Markovian, and hence
so difficult.
 
Whenever a problem seems {\it impossible}, it is not a bad idea to
invent toy problems that
are {\it possible}. Besides the positive chance that it might lead
us eventually to solving the real thing, it is plain fun to play,
and to be able to get results. This is the case for polyominoes.
 
The first toy model for animals was considered by Temperley
([T1],[T2] pp. 66-67)
who treated {\it vertically } (equiv. horizontally) {\it convex}
animals. These are LEGO towers in which every floor consists of
a single piece. Since then, many other, much deeper, toy
models were solved by the \'Ecole Bordelaise (e.g. [DV][V][B1][B2]),
and the Australian school, under the doyenship of Tony
Guttmann (e.g. [BGE]).

In this paper Temperley's method is reviewed, its natural
scope is realized, and then it is generalized in several directions.
Everything has been programmed in Maple, and is contained
in the 
package {\tt LEGO} (available from my homepage), that
can derive Temperley's generating function, and many others,
instantaneously. 
 
{\bf Temperley's Method (slightly rephrased)}
 
Suppose that we have an infinite supply of $1 \times a$ ($a \geq 1$)
LEGO pieces. How many possible {\it LEGO towers} are there, 
with exactly one piece per floor, and
with all the pieces parallel, whose total area is $n$ unit squares?
This problem is equivalent to that of counting {\it vertically convex}
polyominoes, that was solved by Temperley in 1956. His method can
be rephrased as follows. Let $a(n)$ be the required number, then
the generating function $f(t)=\sum_{n=1}^{\infty} a(n)t^n$ is
the {\it weight enumerator} of all compositions (i.e. sequences
of positive integers) with the weight
$$
wt(a_1, \dots , a_r)= t^{a_1+ \dots + a_r}
\prod_{i=1}^{r-1} (a_i+a_{i+1}-1) \quad.
$$
 
Indeed, every horizontally convex polyomino, of height $r$, gives rise to
a composition $(a_1, \dots , a_r )$, where $a_i \geq 1$ is the length of
the $i^{th}$ floor, and to every pair of adjacent floors, of length $a$ and
$b$, there are $a+b-1$ ways of placing the top one on top of the
bottom one.
 
More generally, let $p(a,b)$ be an arbitrary polynomial, and
$L(a)$ an arbitrary affine-linear form $L(a)=c_1 a + c_0$, where
$c_0,c_1$ are integers, $c_1>0$ and $c_0 \geq 0$. 
Define
$$
wt(a_1, \dots , a_r)= t^{L(a_1)+ \dots + L(a_r)}
\prod_{i=1}^{r-1} p(a_i,a_{i+1}) \quad.
$$
It is required to find
$$
f(t):=\sum_{C \in \C} wt(C) \quad .
$$
Here $\C$ is the set of all compositions.
 
Temperley's {\it trick}, that is hereby promoted to {\it method}, is
to consider the weight enumerators $F(a)$ for the subset of $\C$
whose first component is $a$,
$$
F(a):=\sum_{(a,a_2, \dots , a_r)\in \C} wt(a,a_2, \dots , a_r) \quad,
$$
and consider the two-variable generating function
$$
\Phi(z,t):=\sum_{a=1}^{\infty} F(a)z^a \quad .
$$
Once we know $\Phi$ we would also know $f(t)$, since
$f(t)=\Phi(1,t)$.
 
The natural equations for the $F(a)$ are
$$
F(a)=t^{L(a)}+t^{L(a)} \sum_{b=1}^{\infty} p(a,b) F(b) \quad,
\eqno(*)
$$
since every composition $C$, that starts with $a$, is either $(a)$, whose
weight is $t^{L(a)}$, or is of the form $C=(a,C')$, where
$C'$ is a composition on its own right, whose first component
starts with, say, $b$, for some $b \geq 1$, and then
$wt(C)=t^{L(a)} p(a,b) wt(C')$.
 
Now let's expand the polynomial $p(a,b)$ in powers of $b$:
$$
p(a,b)=\sum_{r=0}^{R} p_r(a)b^r \quad ,
$$
and plug it in $(*)$, to get
$$
F(a)=t^{L(a)}+t^{L(a)} \sum_{b=1}^{\infty}  F(b)
\left ( \sum_{r=0}^{R} p_r(a)b^r \right )
=t^{L(a)}+t^{L(a)} \sum_{r=0}^{R} p_r(a) 
\left ( \sum_{b=1}^{\infty}  F(b) b^r \right ) \quad.
$$
Now multiply both sides by $z^a$, and sum over $a \geq 1$, to get
$$
\sum_{a=1}^{\infty} F(a)z^a=
\sum_{a=1}^{\infty} t^{L(a)}z^a +
\sum_{r=0}^{R} \left ( \sum_{a=1}^{\infty} p_r(a)   t^{L(a)} z^a \right ) 
\left ( \sum_{b=1}^{\infty}  F(b) b^r \right ) \quad.
\eqno(**)
$$
 
Let's define,
$$
h(z,t):=\sum_{a=1}^{\infty} t^{L(a)}z^a \quad , \quad {\rm and} \quad
g_r(t,z):=\sum_{a=1}^{\infty} p_r(a)   t^{L(a)} z^a \,\, ,\,\,
( 0 \leq r \leq R),
$$
which are certain explicitly computable rational functions in $(t,z)$.
Also define
$$
\Phi^{[r]}(z,t):=\left ( {z {{d} \over {dz}}}\right )^r  \Phi(z,t) \quad .
$$
Eq. $(**)$ now becomes:
$$
\Phi(z,t)=
h(z,t)+ \sum_{r=0}^{R} g_r(z,t) \Phi^{[r]}(1,t)  \quad .
\eqno(***)
$$
Now apply $(z {{d} \over {dz}})^s$ to both sides of $(***)$, 
for $s=0,1, \dots , R$, and then plug in $z=1$, in order to
get $R+1$ linear equations, with coefficients that are rational
functions of $t$, for the $R+1$ unknowns
$\Phi^{[r]}(1,t)$, $0 \leq r \leq R$, solve them, and get
in particular $\Phi(1,t)=f(t)$. \halmos
 
The procedure that implements this in {\tt LEGO} is
{\tt LEGO}. After you downloaded {\tt LEGO} to
your working directory, get into maple and
type {\tt read LEGO:}.
The function call is
{\tt  LEGO(L,p,a,b,t)}. For example, to get Temperley's
original generating function, type
{\tt  LEGO(a,a+b-1,a,b,t);}, and Maple would respond with
$$
{{ { t(t-1)^3}\over {4\,t^{3}-7\,t^{2}+5\,t-
1}}}\quad .
$$
 
Another simple example (attributed to Moser in [K])
is to the enumeration of leftist (equivalently
rightist) horizontally convex polyominoes. Here the leftmost point
of every floor is farther to the left (or right on top)
than the leftmost point of the floor below it. Here $p=b$ and the
function call is: {\tt  LEGO(a,b,a,b,t);}, producing the output:
${{t(1-t)} \over {1-3t+t^2}}$, which equals (why?)
$\sum_{n=1}^{\infty} F_{2n-1} t^n$,
and hence the number of leftist horizontally convex polyominoes with
area $n$ is $F_{2n-1}$. 
 
To get the number of LEGO towers, with one piece per floor, and
where every floor  is {\it perpendicular} to the floor below it
(so we have a kind of zig-zag pattern, that lives in three
dimensions), do
{\tt LEGO(a,a*b,a,b,t);}, and you would get:
$$
{{t(1-3t+2t^2-t^3)} \over {1-5t+6t^2-3t^3+t^4}} \quad .
$$
 
The reader is welcome to think up other kinds of towers, and to
use {\tt LEGO} to enumerate them.
 
{\bf Building With Colored Pieces}
 
If we have 
$s$ different colors of LEGO pieces, $i=1, \dots , s$, and that
the area of a piece of color $i$ and length $a$ is given by the affine-
linear function $L_i(a)$. Suppose also that the number of ways of
placing a piece of color $j$ and length $b$ on top of a piece of
color $i$ and length $a$ is $p_{i,j}(a,b)$, where the
$p_{i,j}(a,b)$,
($1 \leq i,j \leq s$) are polynomials in $(a,b)$.
 
This amounts to weighted counting of {\it colored compositions}.
A {\it colored composition} is a sequence of colored integers
$(a_1^{(k_1)}, a_2^{(k_2 )}, \dots , a_i^{(k_i )}, \dots , a_r^{(k_r)})$,
where the superscripts denote the colors, (so $1 \leq k_i \leq s$), and
$$
wt(a_1^{(k_1)}, a_2^{(k_2 )}, \dots , a_i^{(k_i )}, \dots , a_r^{(k_r)})=
t^{L_{k_1}(a_1)+ \dots + L_{k_r}(a_r)}
\prod_{i=1}^{r-1} p_{k_i,k_{i+1}}(a_i,a_{i+1}) \quad.
$$
It is required to find
$$
f(t):=\sum_{C \in \C\C_s} wt(C) \quad .
$$
Here $\C\C_s$ is the set of all colored compositions with $s$ colors.
 
For each color $i$ ($i=1, \dots ,s $), and for each integer $a \geq 1$,
consider the weight enumerators $F_i(a)$ for the subset of $\C\C_s$
whose first component is $a^{(i)}$, that is:
$$
F_i(a):=\sum_{(a^{(i)},a_2^{(k_2)}, \dots , 
a_r^{(k_r)})\in \C\C_s} 
wt(a^{(i)},a_2^{(k_2)}, \dots , a_r^{(k_r)}) \quad,
$$
and define the two-variable generating functions
$$
\Phi_i(z,t):=\sum_{a=1}^{\infty} F_i(a)z^a \quad (1 \leq i \leq s)\,\, .
$$
Once we know the $\Phi_i$ we would also know $f(t)$, since
$f(t)=\sum_{i=1}^{s}\Phi_i(1,t)$.
 
The natural equations for the $F_i(a)$ are
$$
F_i(a)=t^{L_i(a)}+t^{L_i(a)} \sum_{b=1}^{\infty} 
\sum_{j=1}^{s} p_{i,j}(a,b) F_j(b) \quad,
\eqno(*)
$$
since every colored 
composition in $\C\C_s$, that starts with $a^{(i)}$, is either $(a^{(i)})$, 
whose weight is $t^{L_i(a)}$, or is of the form $C=(a^{(i)},C')$, where
$C'$ is a composition on its own right, whose first component
starts with, say, $b^{(j)}$, for some $b \geq 1$, and some color
$j$ ($1 \leq j \leq s$), and then
$wt(C)=t^{L_i(a)} p_{i,j}(a,b) wt(C')$.
 
Expand the polynomials $p_{i,j}(a,b)$ in powers of $b$:
$$
p_{i,j}(a,b)=\sum_{r=0}^{R_{i,j}} p_{i,j}^{(r)}(a)b^r \quad ,
$$
and plug it in $(*)$, to get
$$
F_i(a)=t^{L_i(a)}+t^{L_i(a)} \sum_{j=1}^{s}\sum_{b=1}^{\infty}  F_j(b)
\left ( \sum_{r=0}^{R_{i,j}} p_{i,j}^{(r)}(a)b^r \right )
=t^{L_i(a)}+t^{L_i(a)} \sum_{j=1}^{s}\sum_{r=0}^{R_{i,j}} p_{i,j}^{(r)}(a) 
\left ( \sum_{b=1}^{\infty}  F_j(b) b^r \right ) \quad.
$$
Now multiply both sides by $z^a$, and sum over $a \geq 1$, to get
$$
\sum_{a=1}^{\infty} F_i(a)z^a=
\sum_{a=1}^{\infty} t^{L_i(a)}z^a +
\sum_{j=1}^{s}
\sum_{r=0}^{R_{i,j}} \left 
( \sum_{a=1}^{\infty} p_{i,j}^{(r)}(a)   t^{L_i(a)} z^a \right ) 
\left ( \sum_{b=1}^{\infty}  F_j(b) b^r \right ) \quad.
\eqno(**)
$$
 
Let's define,
$$
h_i(z,t):=\sum_{a=1}^{\infty} t^{L_i(a)}z^a \quad , \quad {\rm and} \quad
g_{i,j}^{(r)}(t,z):=\sum_{a=1}^{\infty} p_{i,j}^{(r)}(a)   
t^{L_i(a)} z^a \,\, ,\,\,
( 0 \leq r \leq R_{i,j}, 1 \leq i,j \leq s),
$$
which are certain explicitly computable rational functions in $(t,z)$.
Eq. $(**)$ now becomes:
$$
\Phi_i(z,t)=
h_i(z,t)+ \sum_{j=1}^{s}\sum_{r=0}^{R_{i,j}} 
g_{i,j}^{(r)}(z,t) \Phi_j^{[r]}(1,t)  \quad,\,\, ( 1 \leq i \leq s) \,.
\eqno(***)
$$
Now apply $(z {{d} \over {dz}})^l$ to both sides of $(***)$, 
for $l=0,1, \dots , \max_jR_{i,j}$, and then plug in $z=1$, in order to
get $\sum_{i=1}^{s} (\max_jR_{i,j}+1)$ linear equations, with coefficients 
that are rational functions of $t$, for the unknowns
$\Phi_i^{[l]}(1,t)$, $0 \leq l \leq \max_jR_{i,j}$, solve them, and get
in particular $\Phi_i(1,t)$, for $i=1, \dots ,s$, and finally
$f(t)=\sum_{i=1}^{s} \Phi_i(1,t)$.
 
The procedure in {\tt LEGO} that implements the weighted enumeration of
colored compositions is {\tt muLEGO}, the function call is
{\tt muLEGO(Ls,p,a,b,t);}, where $Ls$ is the list 
$[L_1 , \dots , L_s]$, and $p$ is the list of lists
$$
[[p_{1,1}, p_{1,2} , \dots , p_{1,s}],\dots,
[p_{i,1}, \dots , p_{i,j} , \dots , p_{i,s}],
\dots,
[p_{s,1}, \dots , p_{s,j} , \dots , p_{s,s}]] \quad ,
$$
$a,b$ are the variable names, and $t$ is the variable chosen
for the generating function. For example, to find the
generating function for {\it locally stable horizontally convex} polyominoes,
(i.e. the center of gravity of every floor is in the
interior of the floor below it) do:\hfill\break
{\tt muLEGO([2*a,2*a-1],[[2*a-1,2*a],[2*a-2,2*a-1]],a,b,t);}, getting the
output ${{t(1+t-t^2-t^3)} \over {1-t-3 t^2}}$.
 
Yet another example is the number of LEGO towers, with one piece
per floor, but now you have an infinite supply of both
$1 \times a$ and $2 \times a$ pieces. Furthermore, the towers are
to be constructed in such a way that all the pieces are parallel to
each other (they each have a designated length-side and width-side,
even $1 \times 1, 1 \times 2, 2 \times 1$, and $2 \times 2$ pieces).
The function call is \hfill\break
{\tt muLEGO([a,2*a],[[a+b-1,2*(a+b-1)],[2*(a+b-1),3*(a+b-1)]],a,b,t);}.
I omit the output (do it yourself!). More generally, to find
the number of such towers where you have an unlimited supply
of pieces of the shape:
$1 \times a, 2 \times a , \dots , R \times a$, (where $a \geq 1$),
you may use the built-in function {\tt MIGDAL(R,t);}. Thus
Temperley's original generating function is, in particular,
{\tt MIGDAL(1,t)}.
 
{\bf Higher Dimensional Structures}
 
Suppose that we have an infinite supply of $a \times b$ pieces
$1 \leq a,b$, how many towers can we build of side-surface-area
$n$ (we don't count the area of the base and top)
where each floor has exactly one piece, and all the `lengths' are
parallel to each other? Now we have vector compositions
$([a_1^{(1)},a_1^{(2)}],[a_2^{(1)},a_2^{(2)}],\dots,[a_r^{(1)},a_r^{(2)}])$,
where the weight is
$$
wt([a_1^{(1)},a_1^{(2)}],[a_2^{(1)},a_2^{(2)}],\dots,[a_r^{(1)},a_r^{(2)}])=
t^{2 (\sum_{i=1}^{r} a_i^{(1)}+a_i^{(2)})}
\prod_{i=1}^{r-1} (a^{(1)}_i+a^{(1)}_{i+1}-1) 
(a^{(2)}_i+a^{(2)}_{i+1}-1) \quad.
$$
 
More generally, fixing $m$, and using vector notation
$a=(a^{(1)}, \dots , a^{(m)})$, we have to weight-enumerate
sequences $(a_1 , \dots , a_r)$ with the weight given by
$$
wt(a_1, \dots , a_r)= t^{L(a_1)+ \dots + L(a_r)}
\prod_{i=1}^{r-1} p(a_i,a_{i+1}) \quad,
$$
where $L(a)=L(a^{(1)}, \dots , a^{(m)})$ is affine-linear in its
variables and $p(a,b)$ is a polynomial of $2m$ variables. The
previous analysis goes almost verbatim, and is left to the reader.
 
The procedure in {\tt LEGO} that handles this case is
{\tt LEGOmul}. The function call is {\tt LEGOmul(L,p,a,b,t);}.
Here {\tt L} is the affine-linear form in the $m$ variables {\tt a};
{\tt p} is 
a polynomial in the $2m$ variables {\tt a,b}; {\tt a} and {\tt b} are 
the two lists of $m$ variables, used to describe $L$ and $p$;
and {\tt t} is the 
designated variable of the generating function. For example
to solve the problem described above, type: \hfill\break
{\tt LEGOmul(2*a1+2*a2,(a1+b1-1)*(a2+b2-1),[a1,a2],[b1,b2],t);}.

{\bf Building With Colored Multi-dimensional Pieces}
 
Suppose that we have $s$ different colors, $i=1, \dots ,s$,
where the pieces of color $i$ are $d_i$-dimensional pieces,
$a_1 \times a_2 \times \dots \times a_{d_i}$, with 
the $a_1, a_2, \dots , a_{d_i} \geq 1$. The discussion on colored
one-dimensional pieces goes almost verbatim, only now the
$a_i$, $b_i$, and $z_i$ are multi-variables (with $d_i$ variables). 
 
The function call is: {\tt muLEGOmul(Dims,Ls,pols,a,b,t)},
where {\tt Dims} is the list of dimensions of the colors
$i=1, \dots , s$, {\tt Ls} is the list of affine-linear functions, where
{\tt Ls[i]} depends on {\tt a[1], ... , a[Dims[i]]}; {\tt pols} is
the list of lists of polynomials $p_{i,j}(a_i,b_j)$, where
$a_i$ stands for $(a[1], \dots , a[Dims[i]])$, and 
$b_j$ stands for $(b[1], \dots, b[Dims[j]])$; {\tt a} and {\tt b} are
the letters chosen to express the indexed variable; and {\tt t} is
the designated variable name for the argument of the output.
 
Here are two examples:
 
{\tt
muLEGOmul([2,2],[a[1]+a[2],2*a[1]+2*a[2]],
[[(a[1]+b[1]-1)*(a[2]+b[2]-1), (a[1]+b[1]-1)*(a[2]+b[2]-1)],
[(a[1]+b[1]-1)*(a[2]+b[2]-1), (a[1]+b[1]-1)*(a[2]+b[2]-1)]]
,a,b,t);
}
 
{\tt muLEGOmul([1,2],[a[1]+1,a[1]+a[2]],
[[a[1]+b[1]-1,(a[1]+b[1]-1)*b[2]],[a[2]*(a[1]+b[1]-1),(a[2]+b[2]-1)*
(a[1]+b[1]-1)]]
,a,b,t);}.
 
{\bf Future Directions}
 
The present generalizations of Temperley's method should be
extendible much further, for example to the enumeration of
{\it convex polyominoes}, both according to area and perimeter
(see [B1][B2] and references thereof). 
But, now we no longer get rational functions, and
the natural context would be functional and functional-differential
equations, that might, if in luck (like in [DV]) turn out to be
an algebraic generating function, in which case the future program
should be able to guess it empirically, and then prove it
rigorously.
 
Now the natural equations for the $F(a)$ would be:
$$
F(a)=t^{L(a)}+t^{L(a)} \sum_{b=1}^{a} p(a,b) F(b)
+t^{L(a)} \sum_{b=a+1}^{\infty} q(a,b) F(b) \quad,
$$
where $p(a,b)$ and $q(a,b)$ are different polynomials of
$(a,b)$ or $(q_1^a,q_2^a, \dots, q_1^b, q_2^b, \dots )$.
\eject
{\bf REFERENCES}
 
[BGE] R. Brak, A.J. Guttmann, and I.G. Enting,
{\it Exact solutions of the row-convex polygon perimeter generating
function}, J. Phys. A: Math. Gen., {\bf 23}(1990), 2319-2326.

[B1] M. Bousquet-M\'elou, ``{\it q-\'Enum\'eration de polyominos convexes}'',
Publications de LACIM, UQAM, Montr\'eal, 1991.

[B2] M. Bousquet-M\'elou, ``{\it A method for the enumeration of various
classes of column-convex polygons}'',
Discrete Math. {\bf 154}(1996), 1-25.
 
[DV] M.P. Delest and X.G. Viennot, {\it Algebraic languages and polyominoes
enumeration}, Theor. Comp. Sci. {\bf 34}(1984), 169-206.
 
[G1] H.W. Gould,
{\it  ``Research bibliography of two special number sequences''},
Combinatorial Research Institute, Morgantown, W. Va., 1976, 
[MR 53\#5460].
 
[G2] H.W.  Gould,
{\it  Euler's formula for $n$th differences of powers},
  Amer. Math. Monthly {\bf 85} (1978) 450--467. [MR 58\#256]
 
[G3] H.W.  Gould,   {\it Some formulas for Euler's constant},
Bulletin of Number Theory and Related Topics {\bf 10}(1986), 2-9.
[MR 89g\#11011]
 
[K] D. A. Klarner, {\it Some results concerning polyominoes},
Fibonacci Quart. {\bf 3}(1965), 9-20.

[T1] H.N.V. Temperley, {\it Combinatorial problems suggested by
the statistical mechanics of domains and of rubber-like molecules},
Physical Review {\bf 103} (1956), 1-16.
 
[T2] H.N.V. Temperley, {\it ``Graph Theory and Applications''},
Ellis Horwood, Chichester, 1981.
 
[V] X.G. Viennot, {\it Probl\`emes combinatoire pos\'es par la 
physique statistique, S\'eminaire Bourbaki $n^{o}$} 626,
Asterisque {\bf 121-122}(1985), 225-246.
 
\bye